\documentclass[psamsfonts]{amsart}
\usepackage{amssymb}
\usepackage{graphicx,color}
\usepackage{amssymb,amsfonts,amsthm}
\usepackage[all,arc]{xy}
\usepackage{enumerate}
\usepackage{mathrsfs}

\newtheorem{thm}{Theorem}[section]

\newtheorem{appthm}{Theorem}

\newtheorem{lemma}[thm]{Lemma}
\newtheorem{prop}[thm]{Proposition}
\newtheorem{cor}[thm]{Corollary}

\theoremstyle{definition}
\newtheorem{defn}[thm]{Definition}

\theoremstyle{remark}

\makeatletter
\let\c@equation\c@thm
\makeatother
\numberwithin{equation}{section}

\title[Fully Nonlinear Loewner-Nirenberg Problem]{Borderline Case of Fully Nonlinear Loewner-Nirenberg Problem}

\thanks{The author is supported by
National Natural Science Foundation of China (No. 12571212).}

\author{Zhenan Sui}

\address{Institute for Advanced Study in Mathematics of HIT, Harbin Institute of Technology, Harbin, China}
\email{suizhenan@126.com}

\begin{document}

\begin{abstract}
This note proves the nonexistence of complete conformal metric to the borderline case of fully nonlinear Loewner-Nirenberg Problem, which complements Theorem 1.9 of
Z.~Sui, \emph{On fully nonlinear Loewner--Nirenberg problem of Ricci
curvature}, J.~Funct.\ Anal.\ \textbf{286} (2024).
\end{abstract}

\subjclass[2010]{Primary 53C21; Secondary 35J60}

\maketitle

%\tableofcontents

\section {\large Introduction}

\vspace{4mm}

This note focuses on the borderline case of fully nonlinear Loewner-Nirenberg problem. Let $\Omega \subsetneq \mathbb{R}^n$ be a domain in Euclidean space of dimension $n \geq 3$. Let $g$ be the Euclidean metric. Suppose that $\partial \Omega$ is composed by a finite number of disjoint, non-self-intersecting, smooth, compact submanifolds of $\mathbb{R}^n$. We are interested in when a smooth complete metric
$g_u = u^{\frac{4}{n - 2}} g$ with $u > 0$ satisfying
\begin{equation} \label{eq1-1}
\sigma_k \Big( - g_u^{-1} \text{Ric}_{g_u} \Big) := \sigma_k \Big( \lambda \big( - g_u^{-1} \text{Ric}_{g_u} \big) \Big) = 1 \quad \text{  in  } \Omega
\end{equation}
exists. Here $\operatorname{Ric}_g$ is the Ricci tensor of $g$, $\lambda( A ) = ( \lambda_1, \ldots,
\lambda_n )$ is the eigenvalue vector of the matrix $A$, and
\[ \sigma_k (\lambda) =  \sum\limits_{ 1 \leq i_1 <
\cdots < i_k \leq n} \lambda_{i_1} \cdots \lambda_{i_k} \]
is the $k$th elementary symmetric function defined on
G\r arding's cone
\[ \Gamma_k = \{ \lambda  \in \mathbb{R}^n : \sigma_j ( \lambda ) > 0, \, j = 1,
\ldots, k \}. \]

For convenience, set
\[
 a = \frac{n-2}{2},\qquad p = \frac{n+2}{n-2}.
\]
For a positive $C^2$ function, define
\[
 W[u] = (n-2) D^2 u - \frac{n}{u} d u \otimes d u +
       \bigg( \Delta u + \frac{|Du|^2}{u} \bigg) I.
\]
We can see that equation \eqref{eq1-1} can be written as
\begin{equation} \label{eq1-3}
 \sigma_k^{1/k} \big( W[u] \big) = a u^p \qquad \text{ in } \Omega
\end{equation}
because of the relation
\begin{equation*}
\frac{n - 2}{2} \operatorname{Ric}_{g_u} = - (n - 2) \frac{D^2 u}{u} - \bigg( \frac{\Delta u}{u} + \frac{|D u|^2}{u^2} \bigg) g + \frac{n}{u^2} du
\otimes du,
\end{equation*}
where $\operatorname{Ric}_{g_u}$ is the Ricci tensor of $g_u$,
$D u$, $D^2 u$ and $\Delta u$ are respectively the
gradient, Hessian and Laplace-Beltrami operator of $u$ in the Euclidean
metric. In particular, $\operatorname{Ric}_g = 0$.
To guarantee the completeness of the conformal metric, we require the boundary condition
\begin{equation} \label{eq1-2}
u =  \infty \quad \text{ on }  \partial \Omega.
\end{equation}
A $C^2$ function $u$
is admissible in $\Omega$ if
\begin{equation} \label{admissi}
\lambda \big( W[u] \big)(x) \in \Gamma_k \qquad \text{ for any } x \in \Omega.
\end{equation}
Equation \eqref{eq1-3} is elliptic when $u$ is admissible.
Now we can see that solving fully nonlinear Loewner-Nirenberg problem is equivalent to find smooth positive admissible solution to \eqref{eq1-3}--\eqref{eq1-2}.

In \cite{Sui2024}, Sui proved the following result on the solvability of fully nonlinear Loewner-Nirenberg problem of Ricci curvature (see also Gonz\'alez, Li and Nguyen \cite{Gonzalez-Li-Nguyen} for the Schouten tensor case):

\begin{appthm}
Let $\Omega \subsetneq \mathbb{R}^n$ be a domain. Suppose that the boundary $\partial \Omega$ consists of finitely many disjoint non-self-intersecting smooth compact submanifolds, and the maximal solution $u_{\Omega}$ of \eqref{eq1-3} satisfies $u_{\Omega} \not\equiv 0$. If each component of $\partial \Omega$ has codimension $m$ satisfying
\begin{equation*}
v_m = \big( \underbrace{n - m, \ldots, n - m}_{n - m + 1}, \underbrace{2 - m, \ldots, 2 - m}_{m - 1} \big) \in \Gamma_k,
\end{equation*}
then there exists a smooth complete metric $g_u = u^{\frac{4}{n - 2}} g$ satisfying \eqref{eq1-1}, or equivalently, there exists a smooth positive admissible solution to \eqref{eq1-3}--\eqref{eq1-2}. If some component has codimension $m$ satisfying $v_m \notin \overline{\Gamma}_k$, then there does not exist a complete conformal metric satisfying \eqref{eq1-1}, or equivalently, equation \eqref{eq1-3}--\eqref{eq1-2} is not solvable.
\end{appthm}
In particular, \cite{Sui2024} proved
\begin{appthm}
Let $\Omega \subsetneq \mathbb{R}^n$ be a domain whose boundary $\partial\Omega$ consists of finitely many disjoint non-self-intersecting smooth compact hypersurfaces. Then there exists a smooth complete metric $g_u = u^{\frac{4}{n - 2}} g$ satisfying \eqref{eq1-1}, or equivalently, there exists a smooth positive admissible solution $u$ to \eqref{eq1-3}--\eqref{eq1-2}. Moreover, the conformal factor $u$ has the following growth rate
\begin{equation}  \label{growthrate}
{\rho}^{ \frac{n}{2} - 1}(x) u(x) \rightarrow  \big( (n - 1) (C_n^k)^{\frac{1}{k}} \big)^{\frac{n - 2}{4}} \quad \text{ as }  x \rightarrow \partial\Omega,
\end{equation}
where $\rho (x)$ is the distance of $x$ to $\partial \Omega$.
If $\Omega$ is bounded or $k = 1$, the solution $u$ is unique.
\end{appthm}

In this paper, we consider the borderline case $v_m \in \partial{\Gamma_k}$. We have the following nonexistence result.

\begin{thm} \label{main}
Let $\Sigma$ be a smooth compact embedded submanifold of
$\mathbb{R}^n$, without boundary, of codimension $m$. Suppose that a
punctured tubular neighborhood $0 < \rho < r_*$, where
$\rho = \operatorname{dist}(\,\cdot\,,\Sigma)$, is contained in a domain $\Omega$, while $\Sigma \subset \partial\Omega$.
Assume
\[
 v_m \in \partial\Gamma_k.
\]
Then every positive $C^2$ admissible solution of \eqref{eq1-3} in
$\Omega$ is bounded above in some punctured neighborhood of $\Sigma$.
Consequently the metric $u^{4/(n-2)} g$ is incomplete.
\end{thm}
We note that the assumption $v_m \in \partial\Gamma_k$
implies $2 \leq m \leq n-1$.

The $k = 1$ case is the so called Loewner--Nirenberg problem, which was completely solved by Loewner--Nirenberg \cite{Loewner}, Aviles \cite{Aviles} and V\'eron \cite{Veron}.
The Schouten tensor counterpart was solved by Gonz\'alez, Li and Nguyen \cite{Gonzalez-Li-Nguyen} except the borderline case. An essential difference between Ricci curvature case and Schouten curvature case is: the solution of the former is smooth while the solution of the latter is Lipschitz continuous. There are counterexamples to $C^1$ regularity given in Li and Nguyen \cite{LN21} and in Li, Nguyen and Xiong \cite{LNX23}. On smooth compact Riemannian manifold with boundary, the existence result was given in Guan \cite{Guan08} and Gursky--Streets--Warren \cite{Gursky-Streets-Warren} and the growth rate was given in \cite{Gursky-Streets-Warren}.

In previous work \cite{Sui2024}, the parameters $\alpha \beta > a$ in the barrier function are chosen to dominate the rough upper bound $u \leq C \rho^{-a}$. At the borderline cone vector,
this choice loses the required cone separation near $\zeta = 1$. The novelty of this note is:
Lemma \ref{little} supplies precisely the stronger bound needed to set
$\alpha\beta = a$ in the barrier function so that the comparison is possible.

This paper is organized as follows. Section 2 lists some preliminary facts. Section 3 gives an improved growth rate for the borderline case. This improvement allows us to use the barrier constructed in Section 4 in a comparison argument.

\medskip
\noindent
{\bf Acknowledgements} \quad
The author is supported by National Natural Science Foundation of China (No. 12571212).

\vspace{4mm}

\section{Preliminaries}

\vspace{4mm}

We may need the following elementary facts, which can be found in \cite{Gursky-Streets-Warren, Gonzalez-Li-Nguyen, Sui2024}.

\begin{defn} \label{Def1}
A function $0 < \underline{u} \in C^2(\Omega)$ is a subsolution of \eqref{eq1-3} in $\Omega$ if
\[  \lambda\big( W [\underline{u}] \big) \in \Gamma_k  \text{    and    } \sigma_k^{\frac{1}{k}} \big( W [ \underline{u} ] \big)  \geq  \frac{n - 2}{2} \underline{u}^{\frac{n + 2}{n - 2}} \quad \text{ in }  \Omega.     \]
A function $0 < \overline{u} \in C^2(\Omega)$ is a supersolution of \eqref{eq1-3} in $\Omega$ if
\[ \text{either  } \lambda\big( W [\overline{u}] \big) \notin \Gamma_k  \text{ or }
\sigma_k^{\frac{1}{k}} \big( W [ \overline{u} ] \big)  \leq  \frac{n - 2}{2} \overline{u}^{\frac{n + 2}{n - 2}} \quad \text{ in }  \Omega.     \]
\end{defn}

\begin{prop} \label{prop}
If $u$ is a positive subsolution of \eqref{eq1-3}, so does $c u$ for any constant $0 < c < 1$.  If $u$ is a positive supersolution of \eqref{eq1-3}, so does $c u$ for any constant $c > 1$; if in addition $\lambda\big( W [u] \big) \notin \Gamma_k$ everywhere in $\Omega$, then $c u$ is a supersolution of \eqref{eq1-3} for any $c > 0$.
\end{prop}

\begin{thm} \label{MP}
Let $M$ be a smooth compact manifold with boundary. Suppose that $u$ and $v$ are $C^2$ positive subsolution and supersolution of \eqref{eq1-3} respectively on $M$. If $u \leq v$ on $\partial M$, then $u \leq v$ on $M$.
\end{thm}

Let
\[ B_R(x_0) = \big\{ x \in \mathbb{R}^n  :  |x - x_0| < R \big\}. \]

\begin{prop} \label{subsolution}
For any fixed $s > 0$,
\begin{equation*}
 u (x) =  \Big( 4 (n - 1) (C_n^k)^{\frac{1}{k}} s^2 \Big)^{\frac{n - 2}{4}} \Big( s^2 - {|x - x_0|}^2 \Big)^{1 - \frac{n}{2}} \quad \text{ in  }  B_s (x_0),
\end{equation*}
\begin{equation*}
v (x) = \Big( 4 (n - 1) (C_n^k)^{\frac{1}{k}}  s^2 \Big)^{\frac{n - 2}{4}} \Big( {|x - x_0|}^2 - s^2 \Big)^{1 - \frac{n}{2}} \quad \text{  in  } \mathbb{R}^n\setminus \overline{B_s (x_0)}
\end{equation*}
are admissible solutions of \eqref{eq1-3} which approach to $\infty$ on $\partial B_s (x_0)$.
\end{prop}

\begin{cor}
For any positive admissible solution $u$ of \eqref{eq1-3} in $B_R(x_0)$, we have
\begin{equation} \label{ball}
 u(x) \leq \Big( 4 (n - 1) (C_n^k)^{\frac{1}{k}} R^2 \Big)^{\frac{n - 2}{4}} \Big( R^2 - {|x - x_0|}^2 \Big)^{1 - \frac{n}{2}} \quad\text{  in  } B_R(x_0).
\end{equation}
In particular, we have
\begin{equation} \label{eq1_cor}
 u(x_0) \leq  \Big( 4 (n - 1) (C_n^k)^{\frac{1}{k}} \Big)^{\frac{n - 2}{4}} R^{1 - \frac{n}{2}}.
\end{equation}
\end{cor}

\vspace{4mm}

\section{A vanishing normalized growth rate}

The ball comparison \eqref{eq1_cor} gives, uniformly near
\(\Sigma\),
\begin{equation}\label{rough}
  u(x)\leq C\rho(x)^{-a}.
\end{equation}
One can use a ball of radius \(\rho(x)/2\) centered at \(x\).
The following result gives an improved estimate.

\begin{lemma} \label{little}
Under the hypotheses of Theorem \ref{main},
\[
 \lim_{r\downarrow0}\ \sup_{0<\rho(x)<r}\rho(x)^a u(x)=0.
\]
\end{lemma}

\begin{proof}
By \eqref{rough} the limit on the left exists and is finite. Denote it
by $M$. Suppose $M > 0$. We can choose $x_j$ with
$0 < r_j = \rho(x_j) \to 0$ and $r_j^a u(x_j) \to M$. Let $p_j$ be the
nearest point of $\Sigma$ to $x_j$. Now we choose orthogonal maps $Q_j$ which takes the
fixed plane
\[
 E = \mathbb{R}^{n-m} \times \{0\} \subset \mathbb{R}^{n-m} \times \mathbb{R}^m
\]
to $T_{p_j} \Sigma$, and a fixed unit vector $e \in E^\perp$ to
$(x_j - p_j)/r_j$. Set
\[
 u_j(y) = r_j^a u( p_j + r_j Q_j y ).
\]
Equation \eqref{eq1-3} is invariant under this transformation, that is,
\[ \sigma_k^{1/k} \big( W[u_j] (y) \big) = a u_j^p (y)  \]
because
$ap = a + 2$ and
\[
 W[u_j](y) = r_j^{a + 2} Q_j^T W[u](p_j + r_j Q_j y) Q_j.
\]

For each $j$, define the rescaled submanifold
\[
\Sigma_j := \big\{r_j^{-1}Q_j^{T}(q-p_j): q \in \Sigma\big\}.
\]
Because $Q_j$ is orthogonal and $\rho(x) = \operatorname{dist}(x,\Sigma)$, we have
\[
\begin{aligned}
d_j(y)
&:= r_j^{-1} \rho(p_j + r_j Q_j y) \\
&= r_j^{-1} \inf_{q \in \Sigma} |p_j + r_j Q_j y - q| \\
&= \inf_{q \in \Sigma}
  \big| y - r_j^{-1} Q_j^{T} (q - p_j) \big| \\
&= \operatorname{dist}(y, \Sigma_j).
\end{aligned}
\]
By Lemma \ref{Lemma1} in the appendix, we know that
$d_j(y) = \operatorname{dist}(y, \Sigma_j)$ converges uniformly on compact sets to
$d_E(y) := \operatorname{dist}(y, E)$.

We observe that
\[ u_j(e) = r_j^a u( p_j + r_j Q_j e ) = r_j^a u( x_j ) \to M \quad \text{ as } j \to \infty. \]
In a fixed closed ball $\overline B$ centered at $e$ of radius
less than $1/2$, the definition of $M$ therefore implies the
uniform upper estimate
\begin{equation} \label{upper}
 u_j(y) \leq \big( M + o(1) \big) d_j(y)^{-a}
          = M d_E(y)^{-a} + o(1).
\end{equation}
Indeed, for $y \in \overline B$, put $x = p_j + r_j Q_j y$. Since
$\rho(x) = r_j d_j(y)$, the rescaling gives the identity
\[
 u_j(y) = d_j(y)^{-a} \rho(x)^a u(x).
\]
Also,
\[ \begin{aligned}
\rho(x) = & \rho(p_j + r_j Q_j y) = \rho(p_j + r_j Q_j y) - \rho(p_j) \\
\leq & r_j |Q_j y| = r_j |y| \to 0 \quad \text{ as } j \to \infty \quad \text{ uniformly on } \overline B.
\end{aligned} \]
Thus the definition of $M = \limsup\limits_{\rho \downarrow 0} \rho(x)^a u(x)$
yields $\rho(x)^a u(x) \leq M + o(1)$, uniformly for $y \in \overline B$.
This proves the first inequality in \eqref{upper}.
Also $d_E \geq 1/2$ on $\overline B$. Hence $d_j^{-a} \to d_E^{-a}$ uniformly on $\overline B$ and the second equality in \eqref{upper} is verified.

Define the smooth positive function
\[
 \phi(y) = M d_E(y)^{-a}
\]
on a neighborhood of $\overline B$, and fix $\eta > 0$.
We observe that
\[
 D^2 \big( |y-e|^4 \big) = 4 |y-e|^2 I + 8 (y-e) \otimes (y-e).
\]
Let $y_j$ maximize
\[
 u_j(y) - \phi(y) - \eta |y-e|^4
\]
on $\overline B$, and let its maximum be $c_j$.
By estimate \eqref{upper},
\[ u_j(e) - M = u_j(e) - \phi(e) \leq c_j \leq u_j(y_j) - \phi(y_j) \leq o(1). \]
Hence $c_j \to 0$. Moreover, $y_j \to e$ since
\[ 0 \leq \eta |y_j - e|^4 = u_j (y_j) - \phi(y_j) - c_j \leq o(1) + o(1). \]
In particular the maximum is interior in $B$ for large $j$. The function
\[
 \phi_j(y) := \phi(y) + \eta |y-e|^4 + c_j
\]
touches $u_j$ from above at $y_j$, with equal positive values and
equal gradients. Hence
\[
 W[\phi_j](y_j) - W[u_j](y_j)
 = (n-2) D^2(\phi_j - u_j) + \Delta (\phi_j - u_j) I \geq 0.
\]
It follows that
\[
 \lambda \big( W[\phi_j](y_j) \big) \in \Gamma_k \]
and
\[ \sigma_k^{1/k} \big( W[\phi_j](y_j) \big)
 \geq  \sigma_k^{1/k} \big( W[u_j](y_j) \big)
 =  a u_j(y_j)^p =  a \big(  \phi(y_j) + \eta |y_j - e|^4 + c_j \big)^p.
\]
Letting $j \to \infty$ yields $W[\phi_j](y_j) \to W[\phi](e)$, and thus
\begin{equation}\label{contra}
 \sigma_k^{1/k} \big( W[\phi](e) \big) \geq a M^p > 0.
\end{equation}

On the other hand, write
$y = (z, t) \in \mathbb{R}^{n-m} \times \mathbb{R}^m$, so that $d_E(y) = |t| = r$ and
$\phi(y) = M r^{-a}$. Direct calculation shows that
\begin{align*}
& D_z \phi = 0, \quad D_t \phi  = - a M r^{-a-2} t, \quad  |D \phi|^2 = a^2 M^2 r^{-2a-2}, \\
& D_{zz}^2 \phi = 0, \quad D_{zt}^2 \phi = D_{tz}^2 \phi = 0, \quad D_{tt}^2 \phi  = a M (a+2) r^{-a-4} t \otimes t - a M r^{-a-2} I_m, \\
& \Delta \phi  = -a M (m-2-a) r^{-a-2}.
\end{align*}
Therefore
\[
 \lambda\big( W[\phi](e) \big)
 = a M \big( \underbrace{n-m,\dots,n-m}_{n-m+1},
          \underbrace{2-m,\dots,2-m}_{m-1} \big)
 = a M v_m \in \partial\Gamma_k,
\]
which contradicts \eqref{contra}.
\end{proof}

\vspace{4mm}

\section{Comparison and nonexistence}

\vspace{4mm}

Now we construct the critical barrier.
Put
\[
 A_m := \operatorname{diag} \big( (n - m) I_{n - m + 1}, (2 - m) I_{m - 1} \big).
\]
Let
\[
 \mathcal C_k := \{A \in \operatorname{Sym}(n): \lambda(A) \in \overline \Gamma_k \}.
\]
This is a closed convex cone containing the positive semidefinite
cone. Since $A_m \in \partial\mathcal C_k$, we know that $2 \leq m \leq n - 1$ in view of
\[ A_{1} = (n - 1) I_n \in \mathcal C_k \qquad \text{and} \qquad A_{n} = \operatorname{diag} \big( 0, \underbrace{2 - n, \cdots, 2 - n}_{n - 1} \big) \notin \overline{\mathcal C}_k.
\]
There is a nonzero
supporting matrix, for example,
\begin{equation*}
 T := \operatorname{diag} \big( \theta I_{n - m + 1}, \omega I_{m - 1} \big), \qquad
 \theta, \omega \geq 0, \qquad (n - m + 1) \theta + (m - 1) \omega = 1
\end{equation*}
such that
\[
 \operatorname{tr}(TA) \geq 0 \quad \text{ for any } A \in \mathcal C_k, \qquad \operatorname{tr}(TA_m) = 0.
\]
Direct calculation shows that
\[ \omega = \frac{n - m}{(n - 2)(m - 1)}, \qquad \theta = \frac{m - 2}{(n - 2)(n - m + 1)} \leq \frac{1}{2}. \]

Choose the fixed exponents
\begin{equation*}
 \alpha = \frac{n-2}{4n} \in (0, 1), \qquad
 \beta = 2n.
\end{equation*}
Obviously, $\alpha \beta = a$.
For $0 < c \leq 1$ and $D \geq 1$, consider
\[
 \psi_{c, D} = (c \rho^{-\alpha} + D)^{\beta}, \qquad
 \zeta = \frac{c \rho^{-\alpha}}{c \rho^{-\alpha} + D}.
\]
The above barrier has the same form as in Gonz\'alez--Li--Nguyen \cite{Gonzalez-Li-Nguyen}, but the parameters may be chosen differently.

Choosing an appropriate coordinate system as in \cite{Gonzalez-Li-Nguyen}, that is, for each point $x$ in some punctured neighborhood of $\Sigma$, there is a unique nearest point on $\Sigma$, which we choose to be the origin, and choose a coordinate frame ordered as tangent directions,
the radial normal direction, and the remaining normal directions.
Then at $x$,
\[ d \rho \otimes d \rho = \operatorname{diag} \big( \underbrace{0, \dots, 0}_{n - m}, 1,
\underbrace{0, \dots, 0}_{m - 1}  \big) \]
and
\[ \rho D^2 \rho = \operatorname{diag} \big( \underbrace{0, \dots, 0}_{n - m + 1},
\underbrace{1, \dots, 1}_{m - 1}  \big) + O(\rho), \qquad \rho \Delta \rho = m - 1 + O(\rho). \]
Direct differentiation gives
\begin{equation*}
\begin{aligned}
& \psi_{c, D}^{\frac{1}{\beta} - 1} W [ \psi_{c, D} ] \\
= & a c \rho^{- \alpha - 2} \bigg( - (n - 2) \rho D^2 \rho  + \Big( \alpha (- n  - 2 \beta + 2 ) \zeta + (n - 2) (\alpha + 1) \Big)  d \rho \otimes d \rho \\
& + \Big( \alpha (2 \beta - 1) \zeta + \alpha + 1 - \rho \Delta \rho \Big) I \bigg).
\end{aligned}
\end{equation*}
We thus obtain
\begin{equation} \label{barrier}
\frac{\rho^2}{a \zeta \psi_{c,D}} W[ \psi_{c,D} ]
= \frac{1}{a c} \rho^{\alpha + 2} \psi_{c, D}^{\frac{1}{\beta} - 1} W [ \psi_{c, D} ]
= A_m + B(\zeta) + E,
\end{equation}
where
\[ B(\zeta) = (1 - \zeta)
\left(
  \begin{array}{ccc}
    (\alpha - n + 2) I  & \, & \, \\
    \, & ( n - 1 ) \alpha & \, \\
   \, & \, & ( \alpha - n + 2) I \\
  \end{array}
\right),
\]
and
\[ \|E\| \leq C \rho,
\]
where $C$ is independent of both $c$ and $D$.

For $0 < \rho < 1$, $0 < c \leq 1$, and $D \geq 1$,
\[
 1 - \zeta = \frac{D \rho^\alpha}{c + D \rho^\alpha}
 \geq \frac{\rho^\alpha}{1 + \rho^\alpha}
 \geq \frac{1}{2} \rho^\alpha.
\]
It follows that
\[ \begin{aligned}
\operatorname{tr} \big( T B(\zeta) \big)
= & - (1 - \zeta) \Big( (n - 2 - \alpha)(1 - \theta) - \theta \alpha (n - 1) \Big) \\
\leq &  - \frac{1}{2} \rho^\alpha \cdot \frac{3(n-2)}{8},
\end{aligned} \]
where the last inequality is true because $\theta \leq 1/2$ and $\alpha = (n-2)/(4n)$.
Hence we conclude that
\[
 \operatorname{tr} \Big( T \big( A_m + B(\zeta) + E \big) \Big)
\leq - \frac{3(n-2)}{16} \rho^\alpha + C \rho < 0
\]
whenever $0 < \rho < \rho_0$ for sufficiently small $\rho_0 > 0$
independent of $c, D$.
The supporting property of $T$ now gives
\begin{equation*}
 \lambda \big( W[\psi_{c,D}] \big) \notin \overline \Gamma_k
 \qquad \text{ for } 0 < \rho < \rho_0.
\end{equation*}
Thus every such $\psi_{c, D}$ is a supersolution in the sense of
Definition \ref{Def1}.

Choose $D \geq 1$ so that
\[
 D^\beta \geq \max_{\{\rho = \rho_0\}} u.
\]
For each fixed $c \in (0, 1]$, Lemma \ref{little} gives a
$\delta_c > 0$ such that
\[
 u(x) \leq c^\beta \rho(x)^{-a}
 \leq \big(c \rho(x)^{- \alpha} + D \big)^\beta = \psi_{c, D}(x)
 \qquad \text{ for } 0 < \rho \leq \delta_c.
\]
On the outer boundary $\rho = \rho_0$, the same inequality follows
from the choice of $D$.
Applying the comparison principle Theorem \ref{MP} on
$\{\delta < \rho < \rho_0\}$ for \(0<\delta<\delta_c\), we obtain
\[
 u \leq \psi_{c, D} \qquad \text{ for } 0 < \rho < \rho_0.
\]
Finally let $c \downarrow 0$, keeping $D, \alpha, \beta, \rho_0$ fixed.
We conclude
\[
 u(x) \leq D^\beta \qquad \text{ for } 0 < \rho(x) < \rho_0.
\]
This proves the boundedness assertion of Theorem \ref{main}.
A short normal segment approaching $\Sigma$ then has finite length:
\[
 \int_0^{\rho_0} u( p + t \nu)^{2/(n-2)} dt
 \leq D^{2 \beta/(n-2)} \rho_0 < \infty,
\]
proving incompleteness.

\vspace{4mm}

\section*{Appendix.}

\vspace{4mm}

\begin{lemma} \label{Lemma1}
$d_j(y) = \operatorname{dist}(y, \Sigma_j)$ converges uniformly on compact sets to
$d_E(y) = \operatorname{dist}(y, E)$.
\end{lemma}
\begin{proof}
Since $\Sigma$ is compact and smooth, there are $\delta, \kappa > 0$,
independent of $j$, such that for every $p \in \Sigma$ there is a
smooth map $\psi_p: T_p\Sigma \cap B_\delta(0) \to (T_p\Sigma)^\perp$
such that $\psi_p(0) = 0$, $D\psi_p(0) = 0$, $\big| D^2\psi_p (\xi) \big| \leq \kappa$ and
$\Sigma$ is locally given by $\big\{ p + \xi + \psi_p(\xi) \ : \xi \in T_p\Sigma \cap B_\delta(0) \big\}$.
By Taylor's formula, $|\psi_p(\xi)| \leq \frac{\kappa}{2} |\xi|^2$ whenever $\xi \in T_p\Sigma \cap B_\delta(0)$.

Setting $\psi_j := \psi_{p_j}$, and
for $\eta \in E$ with $|\eta| < \delta/r_j$, setting
\[ h_j(\eta) := r_j^{-1} Q_j^{-1} \psi_j(r_j Q_j \eta), \]
we find that
\[ \big| h_j(\eta) \big| \leq \frac{\kappa}{2} r_j |\eta|^2. \]

Fix a compact set $K \subset \mathbb{R}^n$. Denote $R := \sup_{y \in K} |y|$.

For $j$ sufficiently large such that $2 R r_j < \delta$. If $y \in \Sigma_j \cap B_{2R}(0)$,
then $x = p_j + r_j Q_j y \in \Sigma\cap B_\delta(p_j)$. Writing
$x - p_j = \xi + \zeta$ with $\xi\in T_{p_j}\Sigma$,
$\zeta\in(T_{p_j}\Sigma)^\perp$, we have $|\xi| \leq |x - p_j| < \delta$. Therefore we have $\zeta = \psi_j(\xi)$. Now we write $\xi = r_j Q_j \eta$ where $\eta \in E$. Then we find
\[ y = \eta + h_j(\eta), \qquad
 |\eta| = r_j^{-1} |\xi| \leq |y| \leq 2R. \]
Hence
\[ \Sigma_j \cap B_{2R}(0)
 \subset \{\eta + h_j(\eta):\ \eta \in E,\ |\eta| \leq 2R\} \subset \Sigma_j \]
and
\[ \varepsilon_j := \sup_{|\eta| \leq 2R} \big| h_j(\eta) \big|
 \leq C \kappa R^2 r_j \to 0. \]

Since $p_j \in \Sigma$, we have $0 \in \Sigma_j \cap E$. Also,
\[ T_0 \Sigma_j = Q_j^{-1} T_{p_j} \Sigma = E. \]
For $y \in K$, we have
\[
\operatorname{dist}(y, \Sigma_j) \leq |y| \leq R.
\]
Let $a \in \Sigma_j$ satisfy $\operatorname{dist}(y, \Sigma_j) =  d_j(y) = |y - a|$. Then we have
\[ |a| \leq |a - y| + |y| \leq 2 R. \]
Then we can write
\[ a = \eta + h_j(\eta) \quad \text{ with } \eta \in E \ \text{ and } |\eta| \leq 2R. \]
Likewise, let $b \in E$ satisfy $\operatorname{dist}(y, E) = d_E(y) = |y - b|$. We have
$|b| \leq 2R$.
Now we can see that
\[ d_j(y) = |y - \eta - h_j(\eta)| \leq |y - b - h_j(b)| \leq d_E(y) + \varepsilon_j, \]
and
\[ d_E(y) \leq |y - \eta| = \big| y - a + h_j(\eta) \big| \leq d_j(y) + \varepsilon_j. \]
Therefore
\[
 \sup_{K} |d_j - d_E| \leq \varepsilon_j \to 0,
\]
which proves the claimed uniform convergence on compact sets.
\end{proof}

\end{document}